\documentclass[12pt]{article}

\usepackage[margin=2.54cm]{geometry}
\usepackage{setspace}
\usepackage{amscd,amssymb, amsmath, setspace, graphicx}
\usepackage{url, verbatim}
\usepackage{hyperref}

\newcommand{\macaulay}{\mbox{\sc Macaulay2}}
\newcommand{\phcpack}{\mbox{\sc PHCpack}}
\newcommand{\mixedvol}{\mbox{\sc MixedVol}}

\newcommand{\phcpackmtwo}{\mbox{\tt PHCpack}}
\newcommand{\phc}{\mbox{\tt phc}}

\title{PHCpack in Macaulay2}

\author{Elizabeth Gross, 
 Sonja Petrovi\'c, 
Jan Verschelde
\footnote{EG and JV are with Department of Mathematics, Statistics, and Computer Science, University of Illinois at Chicago,   Chicago, IL\ 60607. SP is with Department of Statistics,       Pennsylvania State University,           University Park, PA\ 16802.}
\thanks{
This material is based upon work partly supported by the National Science Foundation under Grants No.\ 0713018 and 1115777.
The authors are grateful to  Anton Leykin for his generous help with {\tt NAGtypes.m2}. 
Work on this package was carried out while SP was in residence at the University of Illinois at Chicago, and latest contribution was made during the 2012 {\tt Macaulay2} workshop, supported by the US National  Science Foundation.
}}

\begin{document}

\date{\today}

\maketitle

\begin{abstract} 
The \macaulay\ package \phcpackmtwo\ provides an interface to \phcpack, a general-purpose polynomial system solver that uses homotopy continuation.  The main method is a numerical blackbox solver which is implemented for all Laurent systems.  The package also provides a fast mixed volume computation, the ability to filter solutions, homotopy path tracking, and a numerical irreducible decomposition method. 
As the size of many problems in applied algebraic geometry often surpasses the capabilities of symbolic software, this package will be of interest to those working on problems involving large polynomial systems. 
\end{abstract}

\section{Numerical homotopy continuation and \phcpack}


Many problems in applied algebraic geometry require solving, or counting the solutions of, a large polynomial or rational system. 
\phcpackmtwo\ is an interface to the program \phcpack, one of several efficient polynomial system solvers that use numerical homotopy continuation methods~\cite{Li03}. 

The basic idea behind homotopy continuation is simple: to solve a polynomial system $f({\bf x})=0$, one first constructs a system $g({\bf x})=0$ that is easy to solve and then constructs a homotopy, $H({\bf x}(t)) = (1-t)g({\bf x})+tf({\bf x})$, in order to numerically track paths from known solutions of $g$ (with $t=0$) to the solutions of the target system $f$ (with $t=1$). 

Available since release 1.4 of \macaulay~\cite{macaulay2}, this package is motivated by~\cite{Ley11} and uses the data types defined by Leykin in {\tt NAGtypes.m2}. 
The main function of the package allows a \macaulay\ user to solve a system numerically through a blackbox solver, where the creation of the start system and homotopy continuation is done behind the scenes.  The package also provides a fast mixed volume computation and allows the user to filter solutions, to track solution paths explicitly, and to perform numerical irreducible decompositions. 

\medskip
In fact,
the interface \phcpackmtwo\ offers access to most of the functionality of the software \phcpack, which  has been serving as a development platform for many of the algorithms in numerical algebraic geometry~\cite{SVW05}. 
Computations in this paper were done with \phc\ version 2.3.61 (version 1.0 was archived in~\cite{Ver99}). 
Since version 2.3.13, \phcpack\ contains \mixedvol~\cite{GLW05},
and more recently added features are described in~\cite{Ver10}.
Note that \phcpack\ can solve Laurent systems, so the package  includes a method to convert a rational system to a Laurent polynomial system. 
The underlying polyhedral methods perform well on benchmark problems; in many of those, 
 the mixed volume is computed essentially instantaneously. 

Although \phcpack\ is open source, we follow the idea 
of OpenXM~\cite{OpenXM} and require only that the executable \phc\
is available in the execution path of the computer.

\section{Numerical solutions of a polynomial system}
The main function, {\tt solveSystem}, returns solutions of a system 
of polynomial or rational equations. 
 Solutions are returned using data types from {\tt NAGtypes}: a collection of {\tt Points} which are approximations to all complex isolated solutions, or a {\tt WitnessSet} for positive-dimensional components.  
The following system consists of $21$ polynomial equations in $21$ unknowns,
related to a Gaussian cycle conjecture~\cite[\S7.4, page 159]{DSS09}
in algebraic statistics. 
The corresponding variety is zero-dimensional of degree $67$.   
{\small
\begin{verbatim}
Macaulay2, version 1.4
i1 : CC[x_11,x_12,x_16,x_22,x_23,x_33,x_34,x_44,x_45,x_55,x_56,x_66,
      y_13,y_14,y_15,y_24,y_25,y_26,y_35,y_36,y_46];
i2 : system = {x_11+2*x_12+2/3*x_16-1,       x_23*y_13+5/2*x_12+11/2*x_22,
     x_33*y_13+x_34*y_14+11/2*x_23, x_34*y_13+x_44*y_14+x_45*y_15,
     x_45*y_14+x_55*y_15+(45/2)*x_56, x_56*y_15+(22/3)*x_16+(45/2)*x_66,
     82/7*x_12+17/2*x_22+12/5*x_23-1, x_34*y_24+(14/11)*x_23+(12/5)*x_33,
     x_44*y_24+x_45*y_25+(12/5)*x_34, x_45*y_24+x_55*y_25+x_56*y_26,
     x_56*y_25+x_66*y_26+(82/7)*x_16, (12/5)*x_23+(282/5)*x_33+(102/14)*x_34-1,
     x_45*y_35+(282/5)*x_34+(102/14)*x_44, x_55*y_35+x_56*y_36+(102/14)*x_45,
     x_16*y_13+x_56*y_35+x_66*y_36, 10/1*x_34+205/16*x_44+(30/2)*x_45-1,
     x_56*y_46+(205/16)*x_45+(30/2)*x_55, x_16*y_14+x_66*y_46+(305/25)*x_56,
     305/25*x_45+517/7*x_55+(89/3)*x_56-1, x_16*y_15+517/78*x_56+(89/3)*x_66,
     (450/21)*x_16+(89/3)*x_56+(293/19)*x_66-1};
i3 : time solutions = solveSystem system ;
     using temporary files /tmp/M2-5331-2PHCinput and /tmp/M2-5331-2PHCoutput
     -- used 0.056381 seconds
i4 : # solutions
o4 = 67
\end{verbatim}
}
Solutions are returned as a list, each entry being of  type {\tt Point}, which includes diagnostic information such as the condition number and the value of the path-tracking variable $t$. This allows one to decide if a solution is ``good" by using  {\tt peek}, suppressed here in the interest of space. 
Note also that the names of temporary input/output files allow the user to access the details of the entire \phc\ computation, if desired.
{\small
\begin{verbatim}
i5 : solutions_0
o5 = {-3.34446-1.36293*ii, 2.1944+.682742*ii, -.0664982-.0038353*ii, 
-2.90447-1.02757*ii, -.0074284+.306877*ii, .018283-.00390056*ii, 
-.0018297-.0708941*ii, .23468+.062941*ii, -.13257-.0064993*ii, 
.00187754+.000036983*ii, .083551+.0025807*ii, -.00348342+.000364664*ii, 
12.0202-34.4693*ii, 14.2637-.25899*ii, 6.77557+.029139*ii, 
5.38612-.346571*ii, 9.67526+.18591*ii, 8.34346-.39709*ii, 
10.7841-27.2306*ii, 11.3312+.8239*ii, 20.0039+.37216*ii}
o5 : Point
\end{verbatim}
}

The solutions can be further refined as necessary. To best illustrate refinement, consider the same system as above, but where the rational coefficients have been changed to larger rational numbers. Let {\tt newSolutions} be the solutions of this modified system {\tt newSystem}, which can be found in the Appendix (suppressed here in the interest of space). 

Solutions with coordinates below or above a given tolerance can be extracted
by {\tt zeroFilter} and {\tt nonZeroFilter}, respectively. 
In the following example, we ask for solutions whose $12$th coordinate is effectively zero (i.e., smaller than $10^{-19}$). Then, we confirm this by refining the answer to precision $64$; notice that the $12$th coordinate is now on the order of  $10^{-67}$. 
{\small 
\begin{verbatim}
i9 : smallSolution = zeroFilter(newSolutions,11,1.0e-19) 
o9 = {{.0677823, -.386278, .0204925, -1.44743, .982877, -.366596, -.435274,
       .725281, -.422346, .0841728, .0218581, 2.23306e-20, 46.7882, -12.922,
       -70.411, 8.2731, -10.9958, 202.197, -43.8649, 306.199, 198.688}}
i10 : time smallerSolution=refineSolutions(newSystem,smallSolution,64)
     -- used 0.008859 seconds 
o10 = {{.0677823, -.386278, .0204925, -1.44743, .982877, -.366596, -.435274,
       .725281, -.422346, .0841728, .0218581, -1.4308e-67, 46.7882, -12.922,
       -70.411, 8.2731, -10.9958, 202.197, -43.8649, 306.199, 198.688}}
\end{verbatim}
}
Note that when refining solutions, \phc\ also recomputes input coefficients to a higher precision, since rational coefficients may not always have an exact floating-point representation when the precision is limited. 
%
%

\section{Mixed volume}

If the system has as many equations as unknowns, the mixed volume gives an 
upper bound on the number of isolated solutions with nonzero coordinates.
For sufficiently generic coefficients, this bound is sharp.  
The function {\tt mixedVolume} is illustrated below.
{\small
\begin{verbatim}
i11 : time mixedVolume (system )
using temporary files /tmp/M2-4281-6PHCinput and /tmp/M2-4281-6PHCoutput
     -- used 0.011375 seconds
o11 = 75
\end{verbatim}
}
This polyhedral computation is faster than solving the system and provides an upper bound on the number of complex isolated roots in the torus.  
Computing  the degree is much slower (and we note that it takes just as long to verify that the variety is zero-dimensional):
{\small
\begin{verbatim}
i12 : time degree ideal(system)
  -- used 767.432 seconds
o12 = 67
\end{verbatim}
}
While mixed volume counts solutions on the torus, one can also compute the stable mixed volume, which counts solutions with zero components as well, by using optional inputs to the method {\tt mixedVolume}.
\phc\ offers additional functionality and flexibility, not all of which we can illustrate in this short note.  Most interestingly, {\tt mixedVolume} offers an option to use a start system, and creates a polyhedral homotopy from a random start system to the given system. 
The interested reader is referred  to the documentation of the package for more information.

\section{Numerical irreducible decomposition}

Given a list of generators of an ideal $I$, the package can also compute
a {\tt NumericalVariety} with a {\tt WitnessSet} 
for each irreducible component of $V(I)$.
The example below appears in~\cite{DSS09}.
{\small
\begin{verbatim}
i13 : CC[x11,x22,x21,x12,x23,x13,x14,x24];
i14 : system={x11*x22-x21*x12,x12*x23-x22*x13,x13*x24-x23*x14};
i15 : V=numericalIrreducibleDecomposition(system)
writing output to file /tmp/M2-5241-2PHCoutput
calling phc -c < /tmp/M2-5241-3PHCbatch > /tmp/M2-5241-5PHCsession
output of phc -c is in file /tmp/M2-5241-2PHCoutput
... constructing witness sets ... 
preparing input file to /tmp/M2-5241-7PHCinput
preparing batch file to /tmp/M2-5241-9PHCbatch
... calling monodromy breakup ...
session information of phc -f is in /tmp/M2-5241-10PHCsession
output of phc -f is in file /tmp/M2-5241-8PHCoutput
found 3 irreducible factors 
o15 = A variety of dimension 5 with components in
     dim 5:  (dim=5,deg=4) (dim=5,deg=2) (dim=5,deg=2)
o15 : NumericalVariety
\end{verbatim}
}
Witness sets are accessed by dimension:
{\small
\begin{verbatim}
i16 : WitSets=V#5; 
i17 : w=first WitSets;
i18 : w#IsIrreducible
o18 = true
\end{verbatim}
}
In the above example we found three components of dimension five. Let's verify the solutions.
{\small
\begin{verbatim}
i19 : R=QQ[x11,x22,x21,x12,x23,x13,x14,x24];
i20 : system={x11*x22-x21*x12,x12*x23-x22*x13,x13*x24-x23*x14};
i21 : PD=primaryDecomposition(ideal(system));
i22: PD/print
ideal (x23*x14 - x13*x24, x21*x14 - x11*x24, x22*x14 - x12*x24, 
      x12*x23 - x22*x13, x11*x23 - x21*x13, x11*x22 - x21*x12)
ideal (x13, x23, x11*x22 - x21*x12)
ideal (x12, x22, x23*x14 - x13*x24)
\end{verbatim}
}
As we see, the dimension and degree of each component agree with the numerical calculation:
{\small
\begin{verbatim}
i23 : PD/dim
o23 = {5, 5, 5}
i24 : PD/degree 
o24 = {4, 2, 2} 
\end{verbatim}
}


\smallskip



\bibliographystyle{plain}
\small 
\bibliography{PHCinM2}

\section*{Appendix}

The polynomial system used in some of the examples in Section 2. It has the same support as the example system {\tt system}, but the rational coefficients have been changed. It also has $67$ solutions.
{\small
\begin{verbatim}
i6 : newSystem = {(22531/300)*x_11+(821/70)*x_12+(4507/210)*x_16-1,
       x_23*y_13+(22531/300)*x_12+(821/70)*x_22,
       x_33*y_13+x_34*y_14+(821/70)*x_23, x_34*y_13+x_44*y_14+x_45*y_15,
       x_45*y_14+x_55*y_15+(4507/210)*x_56,
       x_56*y_15+(22531/300)*x_16+(4507/210)*x_66,
       (821/70)*x_12+(140953/11025)*x_22+(12325/504)*x_23-1,
       x_34*y_24+(140953/11025)*x_23+(12325/504)*x_33,
       x_44*y_24+x_45*y_25+(12325/504)*x_34, x_45*y_24+x_55*y_25+x_56*y_26,
       x_56*y_25+x_66*y_26+(821/70)*x_16,
       (12325/504)*x_23+(282013/5184)*x_33+(10231/1440)*x_34-1,
       x_45*y_35+(282013/5184)*x_34+(10231/1440)*x_44,
       x_55*y_35+x_56*y_36+(10231/1440)*x_45, x_16*y_13+x_56*y_35+x_66*y_36,
       (10231/1440)*x_34+(205697/16200)*x_44+(30529/2520)*x_45-1,
       x_56*y_46+(205697/16200)*x_45+(30529/2520)*x_55,
       x_16*y_14+x_66*y_46+(30529/2520)*x_56,
       (30529/2520)*x_45+(5175321/78400)*x_55+(897/35)*x_56-1,
       x_16*y_15+(5175321/78400)*x_56+(897/35)*x_66,
       (4507/210)*x_16+(897/35)*x_56+(293581/19600)*x_66-1};
i7 : time newSolutions=phcSolve(newSystem);
using temporary files /tmp/M2-3582-2PHCinput and /tmp/M2-3582-2PHCoutput
     -- used 0.054471 seconds
o7 : List
i8 : # newSolutions
o8 = 67
\end{verbatim}
}
\end{document}